\documentclass[reqno,10pt]{amsart}
\usepackage{color}
\usepackage{graphicx,color,amssymb,latexsym,amsfonts,amsmath}

\renewcommand{\b}[1]{\mbox{\boldmath $#1$}}

\newtheorem{remark}{Remark}

\newtheorem{theorem}{\bf Theorem}[section]

\numberwithin{equation}{section} \theoremstyle{plain}
\theoremstyle{definition}

\newtheorem{example}{Example}[section]

\DeclareMathOperator*{\argmax}{argmax}

\input cyracc.def

\begin{document}
\title[Large Deviations for
Past-Dependent Recursions] {Large Deviations for Past-Dependent
Recursions}
\author{Klebaner F.}
\address{Department of Statistic,
University of Melbourne,\\ Parkville, Victoria, Australia 3052.
current address School of Mathematical Sciences, Monash
University, Clayton, 3058.}
\email{fima.klebaner@sci.monash.edu.au}.
\author{Liptser R.}
\address{Department of Electrical Engineering Systems, Tel Aviv University,
69978 Tel Aviv, Israel} \email{liptser@eng.tau.ac.il}

\phantom{x} \hfill  {\bf Problems of Information Transmission.}

\hfill {\bf 32}, 1996, No. 4, pp. 23--34.

\hfill {\bf Corrected version} March 2006  \maketitle

\begin{abstract}
\noindent The Large Deviation Principle is established for
stochastic models defined by past-dependent non linear recursions
with small noise. In the Markov case we use the result to obtain
an explicit expression for the asymptotics of exit time.
\end{abstract}

\section{\bf Introduction.}
\label{sec-1} The simplest example of a stochastic model defined
by past-dependent recursion with small noise is a linear model
\begin{eqnarray}
X_k^\varepsilon=\sum_{i=1}^ma_iX_{k-i}^\varepsilon+\varepsilon\xi_k
\label{1.1}
\end{eqnarray}
subject to fixed $X_i^\varepsilon=x_i, i=0,1,...,m-1$, where
$\varepsilon$ is a small parameter and $(\xi_k)_{k\ge m}$ is an
i.i.d. sequence of random variables. In the present paper we
consider a general non linear model:
\begin{eqnarray}
X_k^\varepsilon=f(X_{k-1}^\varepsilon,...,X_{k-m}^\varepsilon,
\varepsilon\xi_k), \label{1.2}
\end{eqnarray}
where $f(z_1,...,z_m,y)$ is a continuous function. Note that the
model \eqref{1.2} includes \eqref{1.1} as a special case. For
$m=1$ the model (\ref{1.2}) defines a discrete time Markov
process. When $\varepsilon\to 0$ random variables
$X_k^\varepsilon$ converge to deterministic ones, say, $X_k$ and
$X_k,k\ge 1$ are determined by the recursion
\begin{eqnarray*}
X_k=f(X_{k-1},...,X_{k-m},0)
\end{eqnarray*}
subject to the same initial condition. Furthermore
$(X_k^\varepsilon)_{k\ge m}$ converges to $(X_k)_{k\ge m}$ in the
metric $\rho(x,y)=\sum_{j\ge m}2^{-j}\frac{|x_j-y_j|}{
1+|x_j-y_j|}$. This fact provides the motivation to consider the
large deviation principle (LDP) for  family
$(X_k^\varepsilon)_{k\ge m}$ in the metric space
$(\mathbb{R}^\infty,\rho)$. For Markov case ($m=1$) the LDP was
considered in \cite{5}, \cite{7} and \cite{8}. The choice of the
metric space $(\mathbb{R}^\infty,\rho)$ is a natural one for
obtaining the LDP for the family $(X_k^\varepsilon)_{k\ge m}$.
Recursion (1.2) defines continuous mapping
$(\varepsilon\xi_k)_{k\ge m}\to (X_k^\varepsilon)_{k\ge m}$ in the
metric $\rho$. This  implies that the LDP for
$(X_k^\varepsilon)_{k\ge m}$ follows from the LDP for
$(\varepsilon\xi_k)_{k\ge m}$ by the continuous mapping method of
Freidlin \cite{3} or contraction principle of Varadhan \cite{11}.
Since $(\xi_k)_{k\ge m}$ is an i.i.d. sequence, the LDP for
$(\varepsilon\xi_k)_{k\ge m}$ holds if it holds for the family
$\varepsilon\xi$, where $\xi$ is a copy of $\xi_k$. It should be
mentioned that not only the rate function but also the rate of
speed $q(\varepsilon)$ depends on the distribution of $\xi$.

\medskip
In Section 4, sufficient conditions proving the LDP  for family
$\varepsilon\xi$ are given. Section 3 contains examples for which
rate functions can be explicitly calculated, and in the Markov
case asymptotics for the probability of $\{\max_{1\le k\le
M}|X_k^\varepsilon|\ge 1\}$ is found (Theorem 3.1). Main results
are formulated in Section 2. One of them gives the asymptotics of
the exit time from the interval  $[-1,1]$ for the Markov family
$(X_k^\varepsilon)_{k\ge 1}$.

\section{\bf Main results}
\label{sec-2}

Following Varadhan \cite{11}, family $(X^{\varepsilon})_{k\ge m}$
is said to satisfy the LDP in the metric space
$(\mathbb{R}^\infty,\rho)$ with the rate of speed $q(\varepsilon)$
and the rate function $J(u)$ if

\medskip
{\bf (0)} there exists a function $J=J(\overline{u}),\overline{u}
=(u_1,u_2,...)\in \mathbb{R}^\infty$ which takes values in

\hskip .23in $[0,\infty]$ such that for every $\alpha\ge 0$ the
set $\Phi(\alpha)=\{\overline{u}\in
R^\infty:J(\overline{u})\le\alpha\}$ is

\hskip .23in compact in $(\mathbb{R}^\infty,\rho);$

\medskip
{\bf (1)} For every closed set $F\in (\mathbb{R}^\infty,\rho)$
$$
\varlimsup_{\varepsilon\to 0}q(\varepsilon) \log
\mathsf{P}((X^{\varepsilon}_k)_{k\ge m}\in F) \le
-\inf_{\overline{u}\in F}J(\overline{u});
$$

\medskip
{\bf (2)} For every open set $G\in (\mathbb{R}^\infty,\rho)$
$$
\varliminf_{\varepsilon\to 0}q(\varepsilon) \log
\mathsf{P}((X^{\varepsilon}_k)_{k\ge m}\in G) \ge
-\inf_{\overline{u}\in G}J(\overline{u}).
$$

\medskip
As was mentioned in Introduction, the LDP for
$(X^{\varepsilon}_k)_{k\ge m}$ is implied by the LDP for family
$\varepsilon\xi$ ($\xi$ is a copy of $\xi_k$). Therefore, we begin
with the LDP $\varepsilon\xi$. Henceforth the following conditions
are assumed to be fulfilled:

\textbf{(A.1)}

- $\mathsf{E}\xi=0$

- $\mathsf{E}e^{t\xi}<\infty, \ t\in \mathbb{R}$ ``Cramer's
condition''.

\textbf{(A.2)} With a cumulant function $H(t)=\log
\mathsf{E}e^{t\xi}$ and the Fenchel-Legendre

\hskip .43in transform $ L(v)=\sup_{t\in \mathbb{R}}[tv-H(t)], $
there exist a function $q(\varepsilon)$, decreasing

\hskip .43in to $0$  as $\varepsilon\downarrow 0$, and a
nonnegative function $I(v)=\lim_{\varepsilon\to
0}q(\varepsilon)L({v/\varepsilon}), v\in \mathbb{R}$

\hskip .43in with properties:

\hskip .43in - $I(0)=0$

\hskip .43in - $\lim_{|v|\to \infty}I(v)=\infty$.

\textbf{(A.3)} If $I(v)<\infty$ for some $v$, then
$t_v^\varepsilon=\argmax \big(t\frac{v}{\varepsilon}-H(t)\big)$ is
finite and
$$
\overline{\lim}_{\varepsilon\to 0}{q(\varepsilon)\over\varepsilon}
|t_v^\varepsilon |<\infty \quad\mbox{and}\quad\lim_{\varepsilon\to
0}\varepsilon^2 H''(t_v^\varepsilon)=0.
$$

\medskip
We notice also that the Cramer condition implies $ H'(0)=0 $ and $
H''(t)\ge 0 $ and left continuity of $I(v)$ in a vicinity of
$v_0=\inf\{v>0:I(v)=\infty\}$ (correspondingly, right continuity
for $v_0<0$).

The main results are given

\begin{theorem}\label{theo-2.1}
Assume {\rm \textbf{(A.1) - (A.3)}}. Then{\rm :}

{\rm 1)} the family $\{\varepsilon\xi\}_{\varepsilon\to 0}$ obeys
the LDP with the rate of speed $q(\varepsilon)$ and the function

\hskip .15in $I(v)$ defined in {\rm \textbf{(A.2)}}{\rm ;}

\medskip
{\rm 2)} the family $\{(\varepsilon\xi_k)_{k\ge
1}\}_{\varepsilon\to 0}$ obeys the LDP in the metric space
$(\mathbb{R}^\infty,\rho)$ with

\hskip .15in the rate function $(\overline{v}= (v_1,v_2,...)\in
\mathbb{R}^\infty)$
\begin{eqnarray}
I_\infty(\overline{v})=\sum_{k=1}^\infty I(v_k); \label{2.9}
\end{eqnarray}

\medskip
{\rm 3)} the family $\{(X_k^\varepsilon)_{k\ge
m}\}_{\varepsilon\to 0}$ obeys the LDP in the metric space
$(\mathbb{R}^\infty,\rho)$ with the rate function
$(\overline{u}=(u_m,u_{m+1},...)\in \mathbb{R}^\infty)$
\begin{eqnarray*}
J_\infty(\overline{u})=
\begin{cases}
\sum\limits_{k=m}^\infty \inf\limits_{v_k:
u_k=f(u_{k-1},...,u_{k-m},v_k)}I(v_k), & u_i=x_i,
i=0,...,m-1\\
\infty, & \mbox{otherwise},
\end{cases}
\end{eqnarray*}
where $\inf(\varnothing)=\infty$.
\end{theorem}

\medskip
\noindent
\begin{remark} {\rm Some time $\textbf{(A.1) - (A.3)}$ might be readily verified
if $\xi$ obeys a decomposition $\xi=\xi^{i}+\xi^{ii}$ with
independent random summands satisfying the Cramer condition. If
for $\xi^i$ the Theorem conditions are satisfied (with
$q^i(\varepsilon)$, $I^i(v)$) and
\begin{eqnarray}
\lim_{\varepsilon\to
0}q^{i}(\varepsilon)L^{ii}({v/\varepsilon})=-\infty, \ v\neq 0,
\label{2.11}
\end{eqnarray}
then the theorem statement is valid the rate of speed
$q(\varepsilon)=q^i(\varepsilon)$ and the rate function
$I(v)\equiv I^i(v)$.

Condition (\ref{2.11}) always holds for random variables with a
finite support.}
\end{remark}

\begin{remark} {\rm The requirement for $f(z_1,...,z_m,y)$ to be continuous
can be relaxed if
$$
H_\varepsilon(t,z_1,...,z_m)=\log
\mathsf{E}\exp\Big(tf(z_1,...,z_m, \varepsilon\xi_1)\Big)
$$
is continuous in $z_1,...,z_m$ for every fixed $t$ and
$\varepsilon$, and there exists a norming factor $q(\varepsilon)$,
that is,
$$
\lim_{\varepsilon\to 0}q(\varepsilon)\sup_{t\in
R}[tu-H_\varepsilon( z_1,...,z_m)]=I(u,z_1,...,z_m).
$$
Then,  the LDP for the family $\{(X_k^\varepsilon)_{k\ge
m}\}_{\varepsilon\to 0}$ may hold with
 rate function
$$
J_\infty(u_1,u_2,...)=\sum_{k\ge m}I(u_k,u_{k-1},...,u_{k-m}).
$$
}
\end{remark}

\medskip
The  asymptotics for exit time is our next result. Let
\begin{eqnarray*}
X_k^\varepsilon=aX_{k-1}^\varepsilon+\varepsilon\xi_k
\end{eqnarray*}
subject to  $X_0^\varepsilon=0$ and $\xi_1$ is $(0,1$)-Gaussian
random variable. Denote $\tau^\varepsilon$ exit time from the
interval $[-1,1]$,
\begin{eqnarray*}
\tau^\varepsilon=\min\{k\ge 1:|X_k^\varepsilon|\ge 1\}.
\end{eqnarray*}
\begin{theorem}\label{theo-2.2}
If $|a|<1$, then $ \varlimsup_{\varepsilon\to 0}\varepsilon^2\log
\mathsf{E}\tau^\varepsilon\le \frac{1}{2}(1-a^2). $

\end{theorem}

\begin{remark}
{\rm The statement from 1) related to a discrete time version of
the Freidlin-Wentzell result on of the exit time asymptotics   for
diffusion processes \cite{4} while a corresponding  discrete time
version can be found in Kifer, \cite{6}. Unfortunately, we could
not apply Kifer's result since in \cite{6} $X_k^\varepsilon$ take
values in a compact while in our case $X_k^\varepsilon\in
\mathbb{R}.$ Therefore, repeating some details from \cite{6}, we
give a self-contained proof.}
\end{remark}

\section{\bf Examples and Applications.}
\label{sec-3}

\begin{example}
The rate of speed $q(\varepsilon)=\varepsilon^2$ and the rate
function $I(v)=\frac{1}{2}v^2$ correspond to the family
$\{\varepsilon\xi\}_{\varepsilon\to 0}$ with $(0,1)$-Gaussian
random variable with the cumulant function $H(t)=\frac{t^2}{2}$.
At the same time the pair
$q(\varepsilon)=\frak{\varepsilon}{|\log\varepsilon|}$, $ I(v)=|v|
$ correspond to the Poisson random variable $\xi$ with parameter 1
and the cumulant function $H(t)=e^t+e^{-t}-2$.

It is interesting to note that for $\xi=\xi^i+\xi^{ii}$, where
$\xi^i$ and $\xi^{ii}$ are independent random variables:

- $\xi^i$ is the Gaussian$(0,1)$ random variable,

- $\xi^{ii}$ is the Poisson$(1)$ random variable,

\smallskip
\noindent then the LDP for family
$\{\varepsilon\xi\}_{\varepsilon\to 0}$ holds with $
q(\varepsilon)=\frac{\varepsilon}{|\log
\varepsilon|}\quad\mbox{and}\quad I(v)=|v|. $
\end{example}

\begin{example}
For a linear in $y$ function $f(z_1,...,z_m,y)$, involving in
(\ref{1.2}):
\begin{eqnarray*}
f(z_1,...,z_m,y)=a(z_1,...,z_m)+b(z_1,...,z_m)y,
\end{eqnarray*}
with positive $b(z_1,...,z_m)$, and  Gaussian$(0,1)$ random
variable $\xi_1$ the rate function is defined as:
\begin{eqnarray*}
J_\infty(\overline{u})= \left\{\begin{array}{ll} \sum_{k=m}^\infty
\frac{(u_k-a(u_{k-1},...,u_{k-m})^2}{ b^2(u_{k-1},...,u_{k-m})}, &
u_0=x_0,...,u_{m-1}=x_{m-1}
\\
\infty, & \mbox{otherwise}.
\end{array}
\right.
\end{eqnarray*}
It can be shown, in particular, that the above formula for the
rate function is preserved if $b(z_1,...,z_m)$ equals zero for
some $(z_1,...,z_m)$ provided the convention $0/0=0$.

In the case of the Markov model
$X_k^\varepsilon=a(X_{k-1}^\varepsilon)+ b(X_{k-1}^\varepsilon)
\varepsilon \xi_k$, an analogy to Freidlin-Wentzell's result
\cite{4} for the diffusion $dX_t^\varepsilon=a(X_t^\varepsilon)dt+
\varepsilon b(X_t^\varepsilon)\varepsilon dW_t$ ($W_t$ is Wiener
process) holds.  Namely,
\begin{eqnarray}
J_\infty(\overline{u})= \left\{\begin{array}{ll} \frac{1}{
2}\sum_{k=1}^\infty\frac{[u_k-a(u_{k-1})]^2}{b^2(u_{k-1})}, &
u_0=x_0\nonumber\\
\infty, & \text{otherwise}.
\end{array}
\right.
\end{eqnarray}
\end{example}
\medskip
\noindent {\bf 3.} The next result plays an important role in
proving Theorem 2.2, it also has an independent interest. Notation
$P_{x_0}$ will be used for designating `$X_0^\varepsilon=x_0$'.

\medskip
\noindent We consider the model
$$
X^\varepsilon_k=aX^\epsilon_{k-1}+\varepsilon\xi_k.
$$
\begin{theorem}\label{theo-3.1}
 Let the assumptions of Theorem 2.2 be
fulfilled and $M\ge 1$. Then,
$$
\lim_{\varepsilon\to 0}\varepsilon^2\log
\mathsf{P}_0\big(\max_{1\le k\le M}|X_k^\varepsilon|\ge 1\big)
=-\frac{1}{2\sum_{k=0}^{M-1}a^{2k}}.
$$
In particular, for $|a|<$1,
$$
\lim_{\varepsilon\to 0}\varepsilon^2\log \mathsf{P}_0\big(\max_{1\le
k\le M}|X_k^\varepsilon|\ge 1\big)
=-\frac{1}{2}(1-a^2)\Big[1-a^{2M}\Big].
$$
\end{theorem}
\begin{proof}
The family $\{(X^\varepsilon_k)_{k\ge 1}\}_{\varepsilon\to 0}$ obeys the LDP
with the rate of speed $\varepsilon^2$ and the rate function
\begin{eqnarray}
J_\infty(\overline{u})= \left\{\begin{array}{ll} \frac{1}{
2}\sum_{k=1}^\infty[u_k-au_{k-1}]^2, &
u_0=0\nonumber\\
\infty, & \text{otherwise}.
\end{array}
\right.
\end{eqnarray}
Denote by
$$
\frak{A}_M=\{\overline{u}:\exists k\le M \ \text{with} \ |u_k|\ge
1 \ \text{and} \ u_{k+1}=au_k, \ \forall \ k\ge M  \}.
$$
Obviously, any sequence $(\overline{u}_n)_{n\ge 1}$ from
$\frak{A}_M$ converging in the metric $\rho$ has a limit
$\overline{u}_o\in\frak{A}_M$, that is, $\frak{A}_M$ is a closed
set. Hence, due to the LDP,
\begin{equation*}
\varlimsup_{\varepsilon\to 0}\varepsilon^2\log
\mathsf{P}\Big(\max_{1\le k\le M}|X_k^\varepsilon|\ge 1\Big) \le
-\min_{\overline{u}\in\frak{A}_M}J_\infty(\overline{u}).
\end{equation*}
If $\overline{u}'\in \frak{A}_M$, then
$u'_{k+1}=au'_k$, $\forall \ k\ge M$, so that, $J_\infty(\overline{u}')
=J_M(\overline{u}')$.
If $\overline{u}'\in\frak{A}_M$, then, there exist a number
$\tau'=\inf\{1\le k\le M: |u'_k|\ge 1\}$. Since $J_{\tau'}(\overline{u}'),
\le J_M(\overline{u}')$
we may restrict a minimization procedure up to
\begin{equation}\label{'+'}
\min_{\substack{u'_0=0, |u'_{\tau'}|\ge 1\\
|u'_k|<1, k\le \tau'-1\\ \tau'\le M}}J_{\tau'}(\overline{u}')=
\min_{\substack{u'_0=0, |u'_{\tau'}|= 1\\
|u'_k|<1, k\le \tau'-1\\ \tau'\le M}}J_{\tau'}(\overline{u}')
\end{equation}
The proof of \eqref{'+'} is based on the fact that the lower bound for
$
\min\limits_{\substack{u'_0=0, |u'_{\tau'}|\ge 1\\
|u'_k|<1, k\le \tau'-1\\ \tau'\le M}}J_{\tau'}(\overline{u}')
$
is attainable on
$
\min_{\substack{u'_0=0, |u'_{\tau'}|= 1\\
|u'_k|<1, k\le \tau'-1\\ \tau'\le M}}.
$

To this end, by letting
\begin{equation*}
u'_k=au'_{k-1}+w_k, \ u_0=0, \ k\le \tau',
\end{equation*}
we find that $u'_{\tau'}=\sum_{k=1}^{\tau'}
a^{\tau'-k}w_k$ and by $|u_{\tau'}|\ge 1$ and the
Cauchy-Schwartz inequality $ 1\le
\sqrt{\sum_{k=1}^{\tau'} a^{2(\tau'-k)} \sum_{k=1}^\tau w_k^2} $.
In other words,
\begin{gather*}
\sum_{k=1}^{\tau'} w_k^2 \ge \frac{ 1}{ \sum_{k=1}^{\tau'}
a^{2(\tau'-k)}}\ge \frac{ 1}{ \sum_{k=1}^M a^{2(M-k)}}=\frac{ 1}{
\sum_{k=0}^{M-1} a^{2k}},
\end{gather*}
where the equality is attainable for $w_k\equiv Ka^{M-k}$ with
free parameter $K$ is chosen such that to keep $|u'_{\tau'}|=1$ or
$|u'_{\tau'}|\ge 1$ respectively. Thus, both sides of \eqref{'+'} possesses the same (attainable)
lower bound: $ \frac{1}{2\sum_{k=0}^{M-1} a^{2k}} $ and, due to
the LDP,
\begin{equation*}
\varlimsup_{\varepsilon\to 0}\varepsilon^2\log
\mathsf{P}\Big(\max_{1\le k\le M}|X_k^\varepsilon|\ge 1\Big) \le
-\frac{1}{2\sum_{k=0}^{M-1} a^{2k}}.
\end{equation*}

In order to prove the lower bound
\begin{equation*}
\varliminf_{\varepsilon\to 0}\varepsilon^2\log
\mathsf{P}\Big(\max_{1\le k\le M}|X_k^\varepsilon|\ge 1\Big) \ge
-\frac{1}{2\sum_{k=0}^{M-1} a^{2k}},
\end{equation*}
we introduce an open subset of $\frak{A}_M$:
$$
\frak{A}^o_M=\{\overline{u}:\exists k\le M \ \text{with} \ |u_k|> 1
\ \text{and} \ u_{k+1}=au_k, \ k\ge M  \}.
$$
So, due to the LDP,
\begin{equation*}
\varliminf_{\varepsilon\to 0}\varepsilon^2\log
\mathsf{P}\Big(\max_{1\le k\le M}|X_k^\varepsilon|\ge 1\Big) \ge-\min_{\overline{u}
\in\frak{A}^o_M}J_\infty(\overline{u}).
\end{equation*}
As previously, $\min_{\overline{u}
\in\frak{A}^o_M}J_\infty(\overline{u})=\min_{\overline{u}
\in\frak{A}^o_M}J_M(\overline{u})$ and, moreover,
$$
\min_{\overline{u}
\in\frak{A}^o_M}J_M(\overline{u})=\frac{1}{ 2}\min_{\substack{u_0=0, |u_\tau|>1\\
|u_k|<1, k\le \tau-1\\
\tau\le M}}\sum_{k=1}^\tau(u_k-au_{k-1})^2
=\min_{\substack{u'_0=0, |u'_{\tau'}|= 1\\
|u'_k|<1, k\le \tau'-1\\ \tau'\le M}}J_{\tau'}(\overline{u}').
$$

\medskip
So, it is left to notice that for $|a|<1$,
$$
\sum_{k=0}^{M-1}a^{2k}=\frac{1}{1-a^2}\big[1-a^{2M}\big].
$$
\end{proof}

\section{\bf LDP for $\b{\varepsilon\xi}$}
\label{sec-4}

In this Section, we prove statement 1) of Theorem 2.1.

There are different approaches for proving the LDP (see e.g.
\cite{1}, \cite{2}, \cite{9}). In our case, following Puhalskii's
main theorem, (\cite{P1}, see also   Theorem 1.3 in \cite{10}), it
suffices to prove the exponential tightness:
\begin{eqnarray}
\lim_{c\to \infty}\varlimsup_{\varepsilon\to 0} q(\varepsilon)\log
\mathsf{P}\big(|\varepsilon\xi|\ge c\big)=-\infty, \label{4.1}
\end{eqnarray}
and the local LDP: for every $v\in\mathbb{R}$
\begin{eqnarray*}
\varlimsup_{\delta\to 0}\varlimsup_{\varepsilon\to 0}
q(\varepsilon)\log
\mathsf{P}\big(|\varepsilon\xi-v|\le\delta\big)=\varliminf_{\delta\to
0}\varliminf_{\varepsilon\to 0} q(\varepsilon)\log
\mathsf{P}\big(|\varepsilon\xi-v|\le\delta\big)=-I(v).
\end{eqnarray*}

\medskip
\noindent (\ref{4.1}) is equivalent to
\begin{eqnarray*}
\lim_{c\to \infty}\varlimsup_{\varepsilon\to 0} q(\varepsilon)\log
\mathsf{P}\big(\pm \varepsilon\xi\ge c\big)--\infty.
\end{eqnarray*}
By the Chernoff inequality $ \mathsf{P}(\varepsilon\xi>c)\le
\exp\big(-(tc)/\varepsilon+H(t)\big) $ and due to
$$
\sup_{t>0}[t(c/\varepsilon)-H(t)]= \sup_{t\in
\mathbb{R}}[t(c/\varepsilon)-H(t)]
$$
and \textbf{(A.2)} we find that $ \varlimsup_{\varepsilon\to
0}q(\varepsilon) \log \mathsf{P}(\varepsilon\xi>c)\le
-I(c)\xrightarrow[c\to\infty]{}-\infty. $ The proof of $
\varliminf_{\varepsilon\to 0}q(\varepsilon) \log
\mathsf{P}(-\varepsilon\xi>c)\le
-I(c)\xrightarrow[c\to\infty]{}-\infty $ is similar.

\medskip
The local LDP is proved in two steps:
\begin{eqnarray*}
&& 1)\quad \varlimsup_{\delta\to 0}\varlimsup_{\varepsilon\to 0}
q(\varepsilon)\log
\mathsf{P}\big(|\varepsilon\xi-v|\le\delta\big)\le I(v)
\\
&& 2)\quad \varliminf_{\delta\to 0}\varliminf_{\varepsilon\to 0}
q(\varepsilon)\log
\mathsf{P}\big(|\varepsilon\xi-v|\le\delta\big)\ge -I(v).
\end{eqnarray*}
Set $Z=\exp\big(t\xi-H(t)\big)$. Since $\mathsf{E}Z=1$, for the
proof of 1) we apply an obvious inequality $
EI(|\xi-{u/\varepsilon}|\le{\delta/\varepsilon})Z\le 1 $ which
remains valid with  $Z$ is replaced by its lower bound
$\underline{Z}=\exp\big(-(\delta/\varepsilon)|t|
+t(u/\varepsilon)-H(t)\big)$ on the set
$\{|\xi-{u/\varepsilon}|\le{\delta/\varepsilon}\}$. The latter
provides
\begin{eqnarray*}
q(\varepsilon)\log \mathsf{P} \big(|\xi-{v/\varepsilon}|\le
{\delta/\varepsilon}\big)\le q(\varepsilon)
(\delta/\varepsilon)|t|-q(\varepsilon)[t(v/\varepsilon)- H(t)].
\end{eqnarray*}
and, due to \textbf{(A.2)} and \textbf{(A.3)}, 1) holds.

\medskip
For 2), it suffices to check only the validity of 2)  for $v$ with
$I(v)<\infty$.

Let us denote $P(y)$ the distribution function of $\xi$. Set
$\Lambda_t(y)=\exp\big(ty-H(t)\big)$. Since
$\int_\mathbb{R}\Lambda_t(y)dP(y)=1$ one van introduce new
distribution function $Q_t(y)$ which obeys the following
properties:
\begin{eqnarray}
\int_RydQ_t(y)=H'(t)\quad\mbox{and}\quad
\int_R[y-H'(t)]^2dQ_t(y)=H''(t). \label{4.8}
\end{eqnarray}
Since $I(v)<\infty$, $t_v^\varepsilon$ is a number.  Then, by
taking $t=t^\varepsilon_v$, we find that
\begin{eqnarray*}
P\big(|\varepsilon\xi-v|\le \delta\big)
&=&\int_{|y-(v/\varepsilon)|\le (\delta/\varepsilon)}
\exp\big(-t_v^\varepsilon y+H(t_u^\varepsilon)\big)
dQ_{t_v^\varepsilon}(x)\nonumber\\
&\ge& \exp\big(-|t_v^\varepsilon|(\delta/\varepsilon)-
t_v^\varepsilon v+H(t_v^\varepsilon)\big)\nonumber\\
& &\times\int_{|x-(v/\varepsilon)|\le (\delta/\varepsilon)}
dQ_{t_u^\varepsilon}(x).
\end{eqnarray*}
Hence, 2) holds if, for instance,
$$
\lim_{\varepsilon\to 0} \int_{|x-(v/\varepsilon)|\le
(\delta/\varepsilon)}dQ_{t_v^\varepsilon}(x)=1, \ \delta>0,
$$
while the latter is verified with the help of Chebyshev's
inequality, \eqref{4.8} and \textbf{(A.3)}:
\begin{eqnarray}
\int_{|x-(v/\varepsilon)|\le
(\delta/\varepsilon)}dQ_{t_v^\varepsilon}(x)
&=&1-\int_{|x-(u/\varepsilon)|>
(\delta/\varepsilon)}dQ_{t_v^\varepsilon}(x)
\nonumber\\
&\ge& 1-\frac{\varepsilon^2}{\delta^2}
\int_R(x-u/\varepsilon)^2dQ_{t_v^\varepsilon}(x)\nonumber\\
&=&1-\frac{\varepsilon^2}{\delta^2}H''(t_v^\varepsilon)\nonumber\\
&\to& 1, \ \ \varepsilon\to 0. \nonumber
\end{eqnarray}

\section{\bf LDP for $\b{(\varepsilon\xi_k)_{k\ge 1}}$}
\label{sec-5}

For $n>1$, the LDP for the family $(\varepsilon\xi_k)_{1\le k\le
n}$ in the metric space $(\mathbb{R}^n,\rho^n)$, where for $x,y\in
R^n$ $\rho^n(x,y)=\sum_{k=1}^n|x_k-y_k|$ with the rate of speed
$q(\varepsilon)$ and the rate function
$$
I_n(v^n)=\sum_{k=1}^nI(v_k)
$$
holds due to the vector $(\varepsilon\xi_k)_{1\le k\le n}$ has
i.i.d. entries (see \cite{12}). Next, by Dawson-G\"artner's
theorem (see \cite{13} or \cite{1}), the LDP for family
$(\varepsilon\xi_k)_{k\ge 1}$ holds with the same rate of speed
and the rate function
\begin{eqnarray*}
I(\overline{v})=\sum_{k=1}^\infty I(v_k).
\end{eqnarray*}

\section{\bf LDP for $\b{(X_k^\varepsilon)_{k\ge m}}$} \label{sec-6}

The mapping $(\varepsilon\xi_k)_ {k\ge m}\to
(X_k^\varepsilon)_{k\ge m}$ is continuous in the metric $\rho$.
Therefore, by the contraction principle (continuous mapping
method) (see \cite{3} and  \cite{11}) the family
$(X_k^\varepsilon)_{k\ge m}$ obeys the LDP with the same rate of
speed and the rate function
$$
J_\infty(\overline{u}) =\inf_{(v_k,k\ge m:
u_k=f(u_{k-1},...,u_{k-m},v_k))}I_\infty(\overline{v}),
$$
where $\inf\{\varnothing\}=\infty$ and $I_\infty(\overline{v})$ is
defined in (\ref{2.9}) and
$$
\inf_{(v_k,k\ge m: u_k=f(u_{k-1},...,u_{k-m},v_k)}
I_\infty(\overline{v}) =\sum_{k=m}^\infty\inf_{(v_k,k\ge m:
u_k=f(u_{k-1},...,u_{k-m},v_k))} I(v_k).
$$

\medskip
Remark 1 holds true since by (\ref{2.11}) the random variable
$\xi_1^{ii}$ is exponentially negligible with respect to the
norming factor $q^i(\varepsilon)$: for any $\delta>0$
$$
\lim_{\varepsilon\to 0}q^i(\varepsilon)\log
\mathsf{P}(|\varepsilon\xi_1^{ii}|> \delta)=-\infty.
$$

\section{\bf Asymptotics of exit time.}
\label{sec-7}

In this Section we prove Theorem \ref{theo-2.2}

Let $M$ be an integer. It is clear that $\varepsilon^2\log
\mathsf{E}\tau^\varepsilon$ and $\varepsilon^2\log
\mathsf{E}\frac{\tau^\varepsilon}{M}$ have the same asymptotics as
$\varepsilon\to 0.$

By taking into account $[z]\le z\le [z]+1$, where $[z]$ is integer
part of $z$, write
\begin{eqnarray*}
\mathsf{E}\frac{\tau^\varepsilon}{M}
&\le& \mathsf{E}\Big[\frac{\tau^\varepsilon}{M}\Big]+1\nonumber\\
&=&\sum_{n=1}^\infty
\mathsf{P}\Big(\Big[\frac{\tau^\varepsilon}{M}\Big]\ge n\Big)+1
\nonumber\\
&\le& \sum_{n=1}^\infty \mathsf{P}\Big(\tau^\varepsilon\ge Mn\Big)+1\nonumber\\
&\le& \sum_{n=0}^\infty \mathsf{P}\Big(\tau^\varepsilon>Mn\Big)+1.
\nonumber
\end{eqnarray*}
This upper bound implies the desired statement if
\begin{eqnarray*}
\lim_{M\to\infty}\lim_{\varepsilon\to 0}\varepsilon^2\log
\sum_{n=0}^\infty \mathsf{P}\Big(\tau> Mn\Big)\le\frac{1}{2}(1-a^2).
\end{eqnarray*}
In order to establish the above upper bound, we use and an obvious equality
$$
\mathsf{P}\big(\tau^\varepsilon>Mn\big) =\mathsf{P}\big(\max_{1\le
k\le Mn}|X_k^\varepsilon|<1\big)
$$
and the Markov property of
$(X_k^\varepsilon)_ {k\ge 1}$:
\begin{gather*}
\mathsf{P}\big(\max_{1\le
k\le Mn}|X_k^\varepsilon|<1\big)
\\
=\mathsf{E}\Big\{I(\max_{1\le k\le M(n-1)}|X_k^\varepsilon|<1)
\mathsf{P}_{X^\varepsilon_{M(n-1)}}\Big(\max_{M(n-1)<k\le
Mn}|X_k^\varepsilon|<1 \Big)\Big\}.
\end{gather*}
The time-homogeneity of $X^\varepsilon_k$ implies
\begin{gather*}
I_{\{|X^\varepsilon_{M(n-1)}|<1\}}\mathsf{P}_{X^\varepsilon_{M(n-1)}}
\Big(\max_{M(n-1)<k\le
Mn}|X_k^\varepsilon|<1 \Big)
\\
\le
I_{\{|X^\varepsilon_{M(n-1)}|<1\}}\sup_{|x|<1}\mathsf{P}_x
\Big(\sup_{0<k\le M}|X_k^\varepsilon|<1 \Big)
\\
\le I_{\{|X^\varepsilon_{M(n-1)}|<1\}}\Big(1-\inf_{|x|<1}\mathsf{P}_x
\Big(\sup_{0<k\le M}|X_k^\varepsilon|\ge 1 \Big)
\\
=I_{\{|X^\varepsilon_{M(n-1)}|<1\}}\Big(1-\mathsf{P}_0
\Big(\sup_{0<k\le M}|X_k^\varepsilon|\ge 1 \Big),
\end{gather*}
where the latter equality is provided by zero mean Gaussian distribution
of $X^\varepsilon_k$, $k=1,\ldots,M$.

Hence we obtain  the recurrence relation
$$
\mathsf{P}\big(\tau^\varepsilon>Mn\big) \le
\mathsf{P}\big(\tau^\varepsilon>M(n-1)\big)\Big(1-
\mathsf{P}_0\big(\max_{1\le k\le M}|X_k^\varepsilon|\ge
1\big)\Big)
$$
with $\mathsf{P}\big(\tau^\varepsilon>0\big)=1$. Iterating it, we
find that for any $n\ge 1$,
$$
\mathsf{P}\big(\tau^\varepsilon>Mn\big) \le \Big(1-
\mathsf{P}_0\big(\max_{1\le k\le M}|X_k^\varepsilon|\ge
1\big)\Big)^n.
$$
Thus,
$$
\sum_{n=0}^\infty \mathsf{P}\big(\tau^\varepsilon>Mn\big) \le
\frac{1}{ \mathsf{P}_0\big(\max_{1\le k\le M}|X_k^\varepsilon|\ge 1\big)}
$$
and it is left to apply Theorem \ref{theo-3.1} \qed

\medskip
\noindent {\bf Acknowledgment.} Authors would like to thank
anonymous reviewer for his comments which  lead to correction of
mistakes and misprints and improved the paper.

\noindent {\bf Added March 2006.} The authors are grateful to
Goran Hognas and Brita Ruths for finding a mistake in Theorem 2.2.
of the original paper, which led to a weaker version of that
theorem in the current revision.


\begin{thebibliography}{99}
\bibitem{1} A. Dembo and O. Zeitouni {\em Large Deviations
Techniques and Applications. } Jones and Bartlet. 1993.
\bibitem{2}J.D. Deuschel and D.W. Stroock  {\em Large Deviations.}
Academic Press. Boston. 1989.
\bibitem{3} M.I. Freidlin  ``Action functional for a class of stochastic
processes''. {\em Theory of Probab. Appl.}, {\bf 17}, 1972, pp.
511--515.
\bibitem{4} M.I. Freidlin and A.D. Wentzell A.D. (1984). {\em Random
Perturbation s of Dynamical Systems.} N.Y. Springer. 1984.
\bibitem{5} Y. Kifer ``A discrete time version of the
Wentzell-Freidlin theory''. {\em Ann. Probab.}, {\bf 18}, 1990,
1676--1692.
\bibitem{6} Y. Kifer, Y.  {\em  Perturbations of Dynamical Systems.}
Birkhauser, Boston. 1988.
\bibitem{7} F.C. Klebaner and O. Zeitouni ``The exit problem for a
class of density-dependent branching systems''. {\em Ann. Appl.
Probab.}, {\bf 4}, No 4, 1994, pp. 1188--1205.
\bibitem{8} G.J. Morrow and S. Sawyer ``Large deviation results for
a class of Markov chains arising from population genetics''. {\em
Ann. Probab.}, {\bf 17}, 1989, pp. 1124--1146.
\bibitem{9} A. Shwarz and A. Weiss {\em Large Deviations for Performance
Analysis, communication and Computing: I, II.}  AT\&T Bell
Laboratories. 1994.
\bibitem{10} R.Sh. Liptser and A.A. Pukhalskii ``Limit theorems on
large deviations for semimartingales''. {\em Stochastic and
Stochastics Reports}, {\bf 38}, 1992, pp. 201--249.
\bibitem{11} S.R.S Varadhan  {\em Large Deviations  and
Applications.} SIAM. Philadelphia. 1984.
\bibitem{12} J. Lynch and J. Sethuraman ``Large deviations
for processes with independent increments''. {\em Ann. Prob.},
{\bf 15}, 1987, pp. 610-627.
\bibitem{13} D.A. Dawson and J. G\"artner ``Large deviations from
the McKean-Vlasov limit for weakly interacting diffusions'', {\em
Stochastics}, {\bf 20}, 1987, pp. 247-308.


\bibitem{P1}  Pukhalskii, A.A.  On functional principle of large
deviations \emph{New trends in Probability and Statistics}.
V.Sazonov and Shervashidze (eds.), Vilnius, Lithuania,
VSP/Mokslas, 1991, pp. 198--218.


\end{thebibliography}
\end{document}